\documentclass[article,12pt]{amsart}
\usepackage{amsfonts}
\usepackage{amsmath}
\usepackage{graphicx} 
\usepackage{color} 

\numberwithin{equation}{section}
\theoremstyle{plain}
\newtheorem{thm}{Theorem}[section]
\theoremstyle{remark}
\newtheorem{rem}{Remark}[section]

\newcommand{\cN}{\mathcal{N}}

\newcommand{\veps}{\varepsilon}
\newcommand{\wh}{\widehat}

\newcommand{\ind}{\mbox{1}\kern-.25em \mbox{I}}
\font\calcal=cmsy10 scaled\magstep1
\def\build#1_#2^#3{\mathrel{\mathop{\kern 0pt#1}\limits_{#2}^{#3}}}

\def\videbox{\mathbin{\vbox{\hrule\hbox{\vrule height1ex \kern.5em
\vrule height1ex}\hrule}}}
\newtheorem{lem}{Lemma}[section]
\newtheorem{defi}[lem]{Definition}

\email{Bernard.Bercu@math.u-bordeaux1.fr}
\email{Victor.Vazquez@siu.buap.mx}

\begin{document}
\title[Strong controllability via the Schur complement]{A new concept of strong
controllability via the Schur complement in adaptive tracking}
\author{Bernard Bercu}
\address{ Universit\'e Bordeaux 1, Institut de Math\'ematiques de Bordeaux,
UMR 5251, 351 cours de la lib\'eration, 33405 Talence cedex, France.}
\author{Victor Vazquez}
\address{ Universidad Aut\'onoma de Puebla, Facultad de Ciencias Fisico
Matem\'aticas, Avenida San Claudio y Rio Verde, 72570 Puebla, Mexico.}

\begin{abstract}
We propose a new concept of strong controllability associated with the Schur
complement of a suitable limiting matrix. This concept allows us to extend
the previous results associated with multidimensional ARX models.
On the one hand, we carry out a sharp analysis of the almost sure convergence
for both least squares and weighted least squares algorithms. On the other hand,
we also provide a central limit theorem and a law of iterated logarithm
for these two stochastic algorithms. Our asymptotic
results are illustrated by numerical simulations.
\end{abstract}

\maketitle

\section{Introduction}


Consider the $d$-dimensional autoregressive process with adaptive control of
order $(p,q)$, $\mbox{ARX}_{d}(p,q)$ for short, given for all $n\geq 0$ by 
\begin{equation}
\label{ARX}
A(R)X_{n+1}=B(R)U_{n}+\veps_{n+1}  
\end{equation}
where $R$ stands for the shift-back operator and $X_{n},U_{n}$ and 
$\varepsilon _{n}$ are the system output, input and driven noise,
respectively. The polynomials $A$ and $B$ are given for all $z\in \mathbb{C}$
by 
\begin{eqnarray*}
A(z) &=&I_{d}-A_{1}z-\cdots -A_{p}z^{p}, \\
B(z) &=&I_{d}+B_{1}z+\cdots +B_{q}z^{q},
\end{eqnarray*}
where $A_{i}$ and $B_{j}$ are unknown square matrices of order $d$ and 
$I_{d} $ is the identity matrix. Denote by $\theta $ the unknown parameter of
the model 
\begin{equation*}
\theta ^{t}=(A_{1},\ldots ,A_{p},B_{1},\ldots ,B_{q}).
\end{equation*}
Relation (\ref{ARX}) can be rewritten as 
\begin{equation}
 \label{MOD}
X_{n+1}=\theta ^{t}\Phi _{n}+U_{n}+\veps_{n+1} 
\end{equation}%
where the regression vector $\Phi _{n}=\left( X_{n}^{p},U_{n-1}^{q}\right) ^{t}$ with 
$X_{n}^{p}=(X_{n}^{t},\ldots ,X_{n-p+1}^{t})$ and 
$U_{n}^{q}=(U_{n}^{t},\ldots ,U_{n-q+1}^{t})$. 
In all the sequel, we shall assume the 
$(\varepsilon _{n})$ is a martingale difference sequence adapted to the
filtration $\mathbb{F}=(\mathcal{F}_{n})$ where $\mathcal{F}_{n}$ stands for
the $\sigma $-algebra of the events occurring up to time $n$. We also assume
that, for all $n\geq 0$, $\mathbb{E}[\varepsilon _{n+1}\varepsilon
_{n+1}^{t}|\mathcal{F}_{n}]=\Gamma $ a.s. where $\Gamma $ is a positive
definite deterministic covariance matrix. \newline
\ 

A wide literature concerning the estimation of $\theta$ as well as on the
tracking control is available, \cite{Astrom}, \cite{Bercu1}, \cite{Bercu2}, 
\cite{Caines}, \cite{Chen1}, \cite{Guo1}, \cite{Guo3}, \cite{Lai}. The
purpose of this paper is to establish sharp asymptotic results for
stochastic algorithms associated with the estimation of $\theta$ via the
introduction of a new concept of strong controllability. The strong
controllability is closely related to the almost sure convergence of the
matrix 
\begin{equation*}
S_{n}=\sum_{k=0}^{n} \Phi_{k}\Phi_{k}^{t}.
\end{equation*}
In the particular case $q=0$, it was shown in \cite{Bercu2} that 
\begin{equation*}
\lim_{n\rightarrow \infty} \frac{S_{n}}{n} = L \hspace{0.5cm}\text{a.s.}
\end{equation*}
where $L$ is the block diagonal matrix of order $dp$ given by 
\begin{equation*}
L=\text{diag}\left(\Gamma, \cdots ,\Gamma\right).
\end{equation*}
Under the classical causality assumption, we shall now prove that 
\begin{equation}  \label{CVGARX}
\lim_{n\rightarrow \infty} \frac{S_{n}}{n} = \Lambda 
\hspace{0.5cm}\text{a.s.}
\end{equation}
where $\Lambda$ is the symmetric square matrix of order $\delta=d(p+q)$ 
\begin{equation*}
\Lambda=\left( 
\begin{array}{cc}
L & K^t \\ 
K & H%
\end{array}
\right)
\end{equation*}
and the matrices $H$ and $K$ will be explicitly calculated. It is well-known 
\cite{Horn} that $\det(\Lambda)=\det(L)\det(S)$ where $S=H-KL^{-1}K^t$ is
the Schur complement of $L$ in $\Lambda$. Moreover, as $L$ is positive
definite, $\Lambda$ is positive definite if and only if $S$ is positive
definite. Via our new concept of strong controllability, we shall propose a
suitable assumption under which $S $ is positive definite. It will allow us
to improve the previous results \cite{Bercu2}, \cite{Bercu3}, \cite{Guo1}, 
\cite{Guo2} by showing a central limit theorem (CLT) and a law of iterated
logarithm (LIL) for both the least squares (LS) and the weighted least
squares (WLS) algorithms associated with the estimation of $\theta$. \newline

The paper is organized as follows. Section $\!2$ is devoted to the
introduction of our new concept of strong controllability together with some
linear algebra calculations. Section $\!3$ deals with the parameter
estimation and the adaptive control. In Section $\!4$, we establish
convergence (\ref{CVGARX}) and we deduce a CLT as well as a LIL for both LS and WLS
algorithms. Some numerical simulations are provided in 
Section $\! 5$. A short conclusion is given in Section $\! 6$. 
All technical proofs are postponed in the Appendices.

\section{Strong controllability}

In all the sequel, we shall make use of the well-known causality assumption
on $B$. More precisely, we assume that for all $z\in \mathbb{C}$ with 
$|z|\leq 1$ 
\begin{equation*}
\det(B(z))\neq 0. \leqno (\text{A}_1)
\end{equation*}
In other words, the polynomial $\det(B(z))$ only has zeros with modulus $> 1$. 
Consequently, if $r>1$ is strictly less than the smallest modulus of the
zeros of $\det(B(z))$, then $B(z)$ is invertible in the ball with center
zero and radius $r$ and $B^{-1}(z)$ is a holomorphic function 
(see e.g. \cite{Duflo} page 155). For all $z\in \mathbb{C}$ 
such that $|z|\leq r$, we shall denote 
\begin{equation}  \label{DEFP}
P(z)=B^{-1}(z)(A(z)-I_d)=\sum_{k=1}^\infty P_kz^k.
\end{equation}
All the matrices $P_k$ may be explicitly calculated as functions of the
matrices $A_i$ and $B_j$. For example, $P_1=-A_1$, $P_2=B_1A_1-A_2$, $
P_3=(B_2-B_1^2)A_1+B_1A_2-A_3$. We shall often make use of the square matrix
of order $dq$ given, if $p\geq q$, by 
\begin{equation*}
\Pi=\left( 
\begin{array}{ccccc}
P_p & P_{p+1} & \cdots & P_{p+q-2} & P_{p+q-1} \\ 
P_{p-1} & P_p & P_{p+1} & \cdots & P_{p+q-2} \\ 
\cdots & \cdots & \cdots & \cdots & \cdots \\ 
P_{p-q+2} & \cdots & P_{p-1} & P_p & P_{p+1} \\ 
P_{p-q+1} & P_{p-q+2} & \cdots & P_{p-1} & P_p
\end{array}
\right)
\end{equation*}
while, if $p\leq q$, by 
\begin{equation*}
\Pi=\left( 
\begin{array}{cccccc}
P_p & P_{p+1} & \cdots & \cdots & P_{p+q-2} & P_{p+q-1} \\ 
\cdots & \cdots & \cdots & \cdots & \cdots & \cdots \\ 
P_{1} & P_2 & \cdots & \cdots & P_{q-1} & P_{q} \\ 
0 & P_1 & P_2 & \cdots & P_{q-2} & P_{q-1} \\ 
\cdots & \cdots & \cdots & \cdots & \cdots & \cdots \\ 
0 & \cdots & 0 & P_1 & \cdots & P_p
\end{array}
\right).
\end{equation*}

\begin{defi}
An $\mbox{ARX}_d(p,q)$ is said to be strongly controllable if $B$ is causal
and $\Pi$ is invertible, 
\begin{equation*}
\det(\Pi)\neq 0. \leqno (\text{A}_2)
\end{equation*}
\end{defi}

\begin{rem}
The concept of strong controllability is not really restrictive. For
example, if $p=q=1$, assumption $(\text{A}_2)$ reduces to $\det(A_1)\neq 0$,
if $p=2$, $q=1$ to $\det(A_2-B_1A_1)\neq 0$, 
if $p=1$, $q=2$ to $\det(A_1)\neq 0$, while if $p=q=2$ to 
\begin{equation*}
\det \left( 
\begin{array}{cc}
A_1 & A_2 -B_1 A_1 \\ 
A_2 - B_1 A_1 & A_3 - B_1 A_2 + (B_1^2 -B_2)A_1
\end{array}
\right) \neq 0.
\end{equation*}
\end{rem}

\ \newline
For all $1\leq i \leq q$, denote by $H_i$ be the square matrix of order $d$ 
\begin{equation*}
H_i=\sum_{k=i}^\infty P_k\Gamma P_{k-i+1}^{t}.
\end{equation*}
In addition, let $H$ be the symmetric square matrix of order $dq$ 
\begin{equation}
\label{DEFH}
H=\left( 
\begin{array}{ccccc}
H_1 & H_2 & \cdots & H_{q-1} & H_q \\ 
H_2^t & H_1 & H_2 & \cdots & H_{q-1} \\ 
\cdots & \cdots & \cdots & \cdots & \cdots \\ 
H_{q-1}^t & \cdots & H_2^t & H_1 & H_2 \\ 
H_q^t & H_{q-1}^t & \cdots & H_2^t & H_1
\end{array}
\right).
\end{equation}
For all $1\leq i \leq p$, let $K_i=P_i\Gamma $ and denote by $K$ the
rectangular matrix of dimension $dq\times dp$ given, if $p\geq q$, by 
\begin{equation*}
K=\left( 
\begin{array}{ccccccc}
0 & K_1 & K_2 & \cdots & \cdots & K_{p-2} & K_{p-1} \\ 
0 & 0 & K_1 & \cdots & \cdots & K_{p-3} & K_{p-2} \\ 
\cdots & \cdots & \cdots & \cdots & \cdots & \cdots & \cdots \\ 
0 & \cdots & 0 & K_1 & K_2 & \cdots & K_{p-q+1} \\ 
0 & \cdots & \cdots & 0 & K_1 & \cdots & K_{p-q}
\end{array}
\right)
\end{equation*}
while, if $p\leq q$, by 
\begin{equation*}
K=\left( 
\begin{array}{ccccc}
0 & K_1 & \cdots & K_{p-2} & K_{p-1} \\ 
0 & 0 & K_1 & \cdots & K_{p-2} \\ 
\cdots & \cdots & \cdots & \cdots & \cdots \\ 
0 & \cdots & 0 & 0 & K_{1} \\ 
0 & 0 & \cdots & 0 & 0 \\ 
\cdots & \cdots & \cdots & \cdots & \cdots \\ 
0 & 0 & \cdots & 0 & 0
\end{array}
\right).
\end{equation*}
Finally, let $L$ be the block diagonal matrix of order $dp$ 
\begin{equation}
\label{DEFL}
L=\left( 
\begin{array}{ccccc}
\Gamma & 0 & \cdots & 0 & 0 \\ 
0 & \Gamma & 0 & \cdots & 0 \\ 
\cdots & \cdots & \cdots & \cdots & \cdots \\ 
0 & \cdots & 0 & \Gamma & 0 \\ 
0 & 0 & \cdots & 0 & \Gamma
\end{array}
\right)
\end{equation}
and denote by $\Lambda$ the symmetric square matrix of order 
$\delta=d(p+q)$ 
\begin{equation}  \label{DEFLAMBDA}
\Lambda=\left( 
\begin{array}{cc}
L & K^t \\ 
K & H
\end{array}
\right).
\end{equation}

The following lemma is the keystone of all our asymptotic results.

\begin{lem}
\label{MAINLEMMA} Let $S$ be the Schur complement of $L$ in $\Lambda$ 
\begin{equation}  \label{SCHUR}
S=H-KL^{-1}K^t.
\end{equation}
If $(\text{A}_1)$ and $(\text{A}_2)$ hold, $S$ and $\Lambda$ are invertible
and 
\begin{equation}  \label{INVLAMBDA}
\Lambda^{-1}= \left( 
\begin{array}{cc}
L^{-1} + L^{-1}K^tS^{-1}KL^{-1} & - L^{-1}K^tS^{-1} \\ 
- S^{-1}KL^{-1} & S^{-1}
\end{array}
\right).
\end{equation}
\end{lem}

\begin{proof}
The proof is given in Appendix\,A.
\end{proof}

\section{Estimation and Adaptive control}

First of all, we focus our attention on the estimation of the parameter $
\theta$. We shall make use of the weighted least squares (WLS) algorithm
introduced by Bercu and Duflo \cite{BercuDuflo}, \cite{Bercu1}, which satisfies,
for all $n\geq 0$, 
\begin{equation}  \label{WLS}
\widehat{\theta}_{n+1}=\widehat{\theta}_{n}+a_nS_{n}^{-1}(a)\Phi_{n}
\left(X_{n+1}-U_{n}-\widehat{\theta}_{n}^{\,t}\Phi_{n}\right){\! }^{t}
\end{equation}
where the initial value $\hat{\theta}_{0}$ may be arbitrarily chosen and 
\begin{equation*}
S_{n}(a)=\sum_{k=0}^{n} a_k\Phi_{k}\Phi_{k}^{t}+I_\delta
\end{equation*}
where the identity matrix $I_\delta$ with $\delta=d(p+q)$ is added in order
to avoid useless invertibility assumption. The choice of the weighted
sequence $(a_{n})$ is crucial. If 
\begin{equation*}
a_{n}=1
\end{equation*}
we find again the standard LS algorithm, while if $\gamma \!> \! 0$, 
\begin{equation*}
a_{n}=\Bigl(\frac{1}{\log s_{n}} \Bigr)^{1+\gamma} 
\hspace{1cm}\text{with}\hspace{1cm} 
s_n=\sum_{k=0}^n||\Phi_k||^2,
\end{equation*}
we obtain the WLS algorithm of Bercu and Duflo.

Next, we are concern with the choice of the adaptive control $U_n$. The
crucial role played by $U_n$ is to regulate the dynamic of the process 
$(X_n) $ by forcing $X_n$ to track step by step a predictable reference
trajectory $x_n$. We shall make use of the adaptive tracking control
proposed by Astr$\ddot{\mbox{o}}$m and Wittenmark \cite{Astrom} given, for
all $n \geq 0$, by 
\begin{equation}  \label{CONTROL}
U_n = x_{n+1}-\widehat{\theta}_n^{\,t}\,\Phi_n.
\end{equation}
By substituting (\ref{CONTROL}) into (\ref{MOD}), we obtain the closed-loop
system 
\begin{equation}  \label{CLS}
X_{n+1} - x_{n+1}= \pi_n + \varepsilon_{n+1}
\end{equation}
where the prediction error 
$\pi_n = (\theta - \widehat\theta_n)^{\,t}\Phi_n$. 
In all the sequel, we assume that the reference trajectory $(x_n)$
satisfies 
\begin{equation}  \label{CT}
\sum_{k=1}^{n} \parallel x_{k} \parallel^{2} =o(n) \hspace{1cm} \text{a.s.}
\end{equation}
In addition, we also assume that the driven noise $(\varepsilon_n)$
satisfies the strong law of large numbers which means that if 
\begin{equation}
\Gamma_{n}=\frac{1}{n}\sum_{k=1}^{n}\varepsilon_{k}\varepsilon_{k}^{t},
\end{equation}
then $\Gamma_{n}$ converges a.s. to $\Gamma$. That is the case if, for
example, $(\varepsilon_n)$ is a white noise or if $(\varepsilon_n)$ has a
finite conditional moment of order $>2$. Finally, let $(C_{n})$ be the
average cost matrix sequence defined by 
\begin{equation*}
C_{n}=\frac{1}{n}\sum_{k=1}^{n}(X_{k}-x_{k})(X_{k}-x_{k})^{t}.
\end{equation*}
The tracking is said to be optimal if $C_{n}$ converges a.s. to $\Gamma$.

\section{Main results}

Our first result concerns the almost sure properties of the LS algorithm.

\begin{thm}
\label{ASPLS} Assume that the $\mbox{ARX}_d(p,q)$ model is strongly
controllable and that $(\varepsilon_n)$ has finite conditional moment of
order $>2$. Then, for the LS algorithm, we have 
\begin{equation}  \label{TH11}
\lim_{n\rightarrow \infty} \frac{S_{n}}{n} = \Lambda \hspace{0.5cm}
\text{a.s.}
\end{equation}
where the limiting matrix $\Lambda$ is given by (\ref{DEFLAMBDA}). In
addition, the tracking is optimal 
\begin{equation}  \label{TH12}
\parallel C_{n}-\Gamma_{n} \parallel = \mathcal{O} \left( \frac{\log n}{n}
\right) \hspace{0.5cm}\text{a.s.}
\end{equation}
We can be more precise in (\ref{TH12}) by 
\begin{equation}  \label{TH13}
\lim_{n\rightarrow \infty} \frac{1}{\log n}\sum_{k=1}^{n}
(X_{k}-x_{k}-\varepsilon_{k}) (X_{k}-x_{k}-\varepsilon_{k})^t 
= \delta \Gamma \hspace{0.5cm}\text{a.s.}
\end{equation}
Finally, $\widehat{\theta}_{n}$ converges almost surely to 
$\theta$ 
\begin{equation}  \label{TH14}
\parallel \widehat{\theta}_{n}-\theta \parallel^{2}= \mathcal{O} 
\left(\frac{\log n}{n} \right) \hspace{0.5cm}\text{a.s.}
\end{equation}
\end{thm}

Our second result is related to the almost sure properties of the WLS
algorithm.

\begin{thm}
\label{ASPWLS} Assume that the $\mbox{ARX}_d(p,q)$ model is strongly
controllable. In addition, suppose that either $(\varepsilon_n)$ is a white
noise or $(\varepsilon_n)$ has finite conditional moment of order $>2$.
Then, for the WLS algorithm, we have 
\begin{equation}  \label{TH21}
\lim_{n\rightarrow \infty} (\log n)^{1+\gamma} \frac{S_{n}(a)}{n} = 
\Lambda \hspace{0.5cm}\text{a.s.}
\end{equation}
where the limiting matrix $\Lambda$ is given by (\ref{DEFLAMBDA}). In
addition, the tracking is optimal 
\begin{equation}  \label{TH22}
\parallel C_{n}-\Gamma_{n} \parallel= o\left( \frac{(\log n)^{1+\gamma}}{n}
\right) \hspace{0.5cm}\text{a.s.}
\end{equation}
Finally, $\widehat{\theta}_{n}$ converges almost surely to  
$\theta$ 
\begin{equation}  \label{TH23}
\parallel \widehat{\theta}_{n}-\theta \parallel^{2}= \mathcal{O} 
\left( \frac{(\log n)^{1+\gamma}}{n} \right) \hspace{0.5cm}\text{a.s.}
\end{equation}
\end{thm}

\begin{proof}
The proof is given in Appendix\,B.
\end{proof}

\begin{rem}
One can observe that Theorems \ref{ASPLS} and \ref{ASPWLS} extend the
results of Bercu \cite{Bercu2} and Guo \cite{Guo2} 
previously established in the AR framework.
\end{rem}

\begin{thm}
\label{CLTLIL} Assume that the $\mbox{ARX}_d(p,q)$ model is strongly
controllable and that $(\varepsilon_n)$ has finite conditional moment of
order $\alpha>2$. In addition, suppose that $(x_n)$ has the same regularity
in norm as $(\varepsilon_n)$ which means that for all $2< \beta < \alpha$ 
\begin{equation}  \label{REGNORM}
\sum_{k=1}^{n}\parallel x_{k}\parallel^{\beta}=\mathcal{O}(n) \hspace{0.5cm}
\text{a.s.}
\end{equation}
Then, the LS and WLS algorithms share the same central limit theorem 
\begin{equation}  \label{TH31}
\sqrt{n} (\widehat{\theta}_{n}-\theta ) 
\build{\longrightarrow}_{}^{{\mbox{\calcal L}}} 
\mathcal{N}(0,\Lambda^{-1}\otimes\Gamma)
\end{equation}
where the inverse matrix $\Lambda^{-1}$ is given by (\ref{INVLAMBDA}) and
the symbol $\otimes$ stands for the matrix Kronecker product. In addition,
for any vectors $u \in \mathbb{R}^{d}$ and $v \in \mathbb{R}^{\delta}$, they
also share the same law of iterated logarithm 
\begin{eqnarray}  \label{TH32}
\limsup_{n \rightarrow \infty} \left(\frac{n}{2 \log\log n} \right)^{1/2}
\!\!v^{t}(\widehat{\theta}_{n}-\theta)u &=& - \liminf_{n \rightarrow \infty}
\left(\frac{n}{2 \log\log n}\right)^{1/2} \!\!v^{t}(\widehat{\theta}
_{n}-\theta)u  \notag \\
&=& \Bigl(v^{t}\Lambda^{-1}v \Bigr)^{1/2}\Bigl(u^{t}\Gamma u \Bigr)^{1/2} 
\hspace{0.5cm}\text{a.s.}
\end{eqnarray}
In particular, 
\begin{equation}  \label{TH33}
\left(\frac{\lambda_{min}\Gamma} {\lambda_{max}\Lambda} \right)\leq
\limsup_{n \rightarrow \infty} \ \left(\frac{n}{2 \log\log n} \right)
\parallel \hat{\theta}_{n}-\theta \parallel^{2} \leq \left(\frac{
\lambda_{max}\Gamma} {\lambda_{min}\Lambda}\right) \hspace{0.5cm}\text{a.s.}
\end{equation}
where $\lambda_{min} \Gamma$ and $\lambda_{max} \Gamma$ are the minimum and
the maximum eigenvalues of $\Gamma$.
\end{thm}

\begin{proof}
The proof is given in Appendix\,C.
\end{proof}

\section{Numerical simulations}

The goal of this section is to propose some numerical experiments
for illustrating the asymptotic results of Section 4.
In order to keep this section brief, we shall only focus our attention 
on a strongly controllable $ARX_d(p,q)$ model
in dimension $d=2$ with $p=1$ and $q=1$. Our numerical simulations 
are based on $M=500$ realizations of sample size $N=1000$. 
For the sake of simplicity, the reference trajectory $(x_n)$
is chosen to be identically zero and the driven noise $(\veps_n)$ is a Gaussian $\cN(0,1)$ white noise.
Consider the $ARX_{2}(1,1)$ model
\begin{equation*}
X_{n+1}=A X_n+U_{n}+B U_{n-1} + \veps _{n+1}
\end{equation*}
where
\begin{equation*}
A =\left( 
\begin{array}{cc}
2 & 0 \\ 
0 & 1
\end{array}
\right) 
\hspace{1cm} \text{and} \hspace{1cm}
B =\frac{1}{4}\left( 
\begin{array}{cc}
3 & \,0 \\ 
0 & \!-2
\end{array}
\right).
\end{equation*}
Figure 1 shows 
the almost sure convergence 
of the LS estimator $\wh{\theta}_n$ to the four coordinates of 
$\theta$ which are different from zero. 
One can observe that $\wh{\theta}_n$ performs very well in the estimation
of $\theta$.
\begin{figure}[!ht] 
\begin{center}
\includegraphics[width=15cm, height=5.5cm]{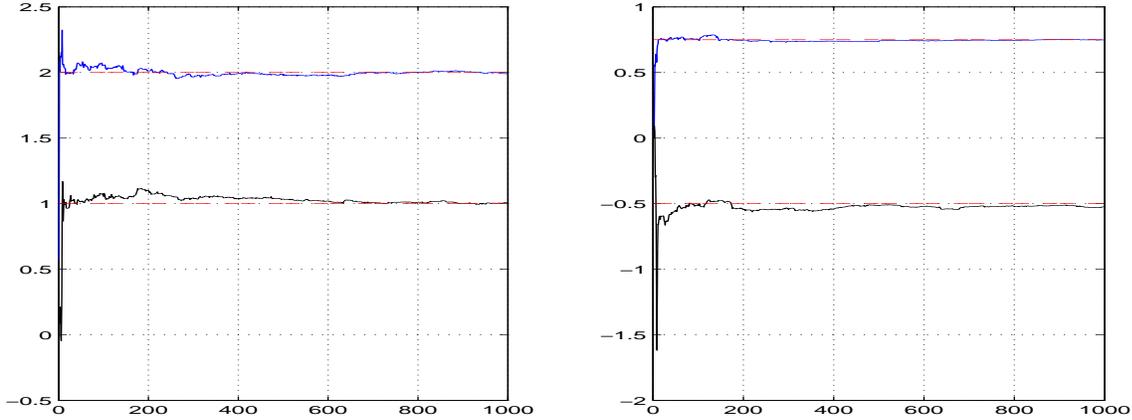} 
\end{center}
\caption{Left: LLN for A, Right: LLN for B.}
\end{figure} 
\ \vspace{1ex}\\
Figure 2 shows the CLT for the four coordinates of 
$$Z_N=\sqrt{N}\Lambda^{1/2}(\wh{\theta}_N - \theta)$$
where $\Lambda$ is the limiting matrix given by (\ref{DEFLAMBDA})
with $L=I_2$, $K=0$ and 
$$H=A^2(I_2-B^2)^{-1}=\frac{4}{21}\left( 
\begin{array}{cc}
48 & 0 \\ 
0 & 7
\end{array}
\right).
$$ 
One can realize that each component
of $Z_N$ has $\cN(0,1)$ distribution as expected.
\begin{figure}[!ht] 
\begin{center}
\includegraphics[width=15cm, height=5.5cm]{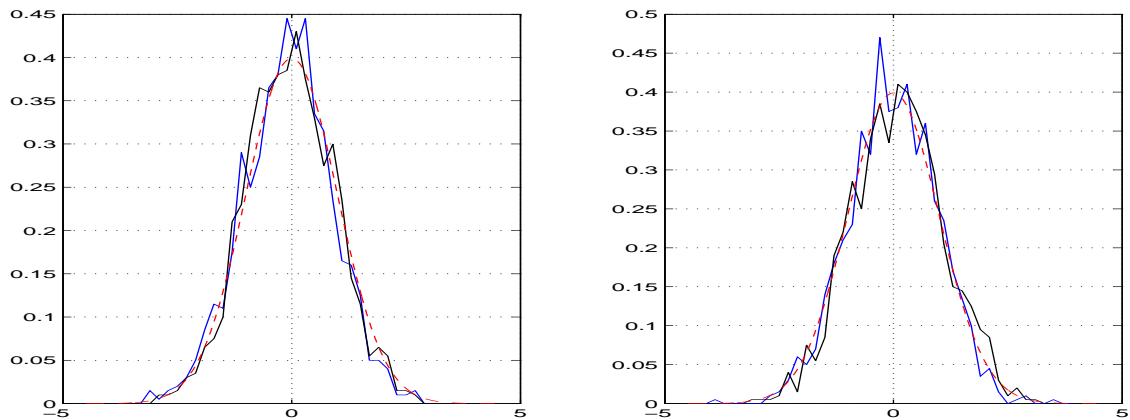} 
\end{center}
\caption{Left: CLT for A, Right: CLT for B.}
\end{figure} 

\section{Conclusion}

Via our new concept of strong controllability, we have extended the analysis of
the almost sure convergence for both LS and WLS algorithms in  
the multidimensional ARX framework. It enables us to provide a positive
answer to a conjecture in \cite{Bercu2} by establishing a CLT and a LIL for these
two stochastic algorithms. In our approach, the leading matrix associated with 
the matrix polynomial $B$, commonly called the high frequency gain, was
supposed to be known and it was chosen as the identity matrix $I_d$.
It would be a great challenge for the control community to carry out 
similar analysis with unknown high frequency gain and to extend it to
ARMAX models.

\section{Appendix A.}

\renewcommand{\thesection}{\Alph{section}} \renewcommand{\theequation}
{\thesection.\arabic{equation}} \setcounter{section}{1} 
\setcounter{equation}{0}


\textit{Proof of Lemma \ref{MAINLEMMA}}. Let $\Sigma$ be the
infinite-dimensional diagonal square matrix 
\begin{equation*}
\Sigma=\left( 
\begin{array}{ccccc}
\Gamma & 0 & \cdots & \cdots & \cdots \\ 
0 & \Gamma & 0 & \cdots & \cdots \\ 
\cdots & 0 & \Gamma & 0 & \cdots \\ 
\cdots & \cdots & \cdots & \cdots & \cdots \\ 
\cdots & \cdots & \cdots & \cdots & \cdots%
\end{array}
\right).
\end{equation*}
Moreover, denote by $T$ the infinite-dimensional rectangular matrix with $dq$
raws and an infinite number of columns given, if $p\geq q$, by 
\begin{equation*}
T=\left( 
\begin{array}{cccccc}
P_p & P_{p+1} & \cdots & P_{k} & P_{k+1} & \cdots \\ 
P_{p-1} & P_p & \cdots & P_{k-1} & P_{k} & \cdots \\ 
\cdots & \cdots & \cdots & \cdots & \cdots & \cdots \\ 
P_{p-q+2} & P_{p-q+3} & \cdots & P_{k-q+2} & P_{k-q+3} & \cdots \\ 
P_{p-q+1} & P_{p-q+2} & \cdots & P_{k-q+1} & P_{k-q+2} & \cdots
\end{array}
\right)
\end{equation*}
while, if $p\leq q$, by 
\begin{equation*}
T=\left( 
\begin{array}{ccccccc}
P_p & P_{p+1} & \cdots & \cdots & P_{k} & P_{k+1} & \cdots \\ 
\cdots & \cdots & \cdots & \cdots & \cdots & \cdots & \cdots \\ 
P_{1} & P_2 & \cdots & \cdots & P_{k-p+1} & P_{k-p+2} & \cdots \\ 
0 & P_1 & P_2 & \cdots & P_{k-p} & P_{k-p+1} & \cdots \\ 
\cdots & \cdots & \cdots & \cdots & \cdots & \cdots & \cdots \\ 
0 & \cdots & \cdots & 0 & P_1 & P_2 & \cdots
\end{array}
\right).
\end{equation*}
After some straightforward, although rather lengthy, linear algebra
calculations, it is possible to deduce from (\ref{SCHUR}) that 
\begin{equation}  \label{SCHURDEC}
S=T \Sigma T^t.
\end{equation}
It clearly follows from this suitable decomposition that 
\begin{equation}  \label{SCHURKER}
\ker(S)=\ker(T^t).
\end{equation}
As a matter of fact, assume that $v \in \mathbb{R}^{dq}$ belongs to 
$\ker(T^t)$. Then, $T^tv=0$, $Sv=0$ which leads to 
$\ker(T^t)\subset \ker(S)$. 
On the other hand, assume that $v \in \mathbb{R}^{dq}$ belongs 
to $\ker(S)$. Since $Sv=0$, we clearly have 
$v^tSv=0$, so $v^t T \Sigma T^t v=0$. However,
the matrix $\Gamma$ is positive definite. Consequently, 
$T^tv=0$ and $\ker(S)\subset \ker(T^t)$, which 
implies (\ref{SCHURKER}). Moreover, it
follows from the well-known rank theorem that 
\begin{equation}  \label{RANKTHM}
dq= \dim(\ker(S))+\text{rank}(S).
\end{equation}
As soon as $\ker(S)=\{0\}$, $\dim(\ker(S))=0$ and we obtain from 
(\ref{RANKTHM}) that $S$ is of full rank $dq$ which means 
that $S$ is invertible.
Furthermore, the left hand side square matrix of order $dq$ of the
infinite-dimensional matrix $T$ is precisely $\Pi$. Consequently, if $\Pi$
is invertible, $\Pi$ is of full rank $dq$, $\ker(\Pi)=\ker(\Pi^t)=\{0\}$ and
the left null space of $T$ reduces to the null vector of $\mathbb{R}^{dq}$.
Hence, if $\Pi$ is invertible, we deduce from (\ref{SCHURKER}) together with
(\ref{RANKTHM}) that $S$ is also invertible. Finally, as 
\begin{equation}  \label{DECDET}
\det(\Lambda)=\det(L)\det(S)=\det(\Gamma)^p\det(S),
\end{equation}
we obtain from (\ref{DECDET}) that $\Lambda$ is invertible and formula 
(\ref{INVLAMBDA}) follows from \cite{Horn} page 18, which completes 
the proof of Lemma \ref{MAINLEMMA}.$\ \ \videbox$


\section*{Appendix B.}

\renewcommand{\thesection}{\Alph{section}} \renewcommand{\theequation}
{\thesection.\arabic{equation}} \setcounter{section}{2} 
\setcounter{equation}{0}

\textit{Proof of Theorem \ref{ASPLS}}. In order to prove 
Theorem \ref{ASPLS}, we shall make use of the same approach than Bercu
\cite{Bercu2} or Guo and Chen \cite{Guo1}. First of all, we recall that
for all $n\geq 0$
\begin{equation}  
\label{EQBASIS}
X_{n+1} - x_{n+1}= \pi_n + \varepsilon_{n+1}.
\end{equation}
In addition, let
\begin{equation*}
s_n=\sum_{k=0}^n||\Phi_k||^2.
\end{equation*}
It follows from (\ref{EQBASIS}) together with 
the strong law of large numbers for
martingales (see e.g. Corollary 1.3.25 of \cite{Duflo}) 
that $n=\mathcal{O}(s_n)$ a.s. Moreover, by Theorem 1 of
\cite{Bercu0} or Lemma 1 of \cite{Guo1}, we have
\begin{equation}  
\label{SUMPIF}
\sum_{k=1}^n(1-f_k)\parallel \pi_k \parallel^2=\mathcal{O}(\log s_n)
 \hspace{0.5cm} \text{a.s.}
\end{equation}
where $f_n=\Phi_n^tS_n^{-1}\Phi_n$. Hence, if $(\varepsilon_n)$ has 
finite conditional moment of order $\alpha>2$, we can show by the causality 
assumption $(\text{A}_1)$ on the matrix polynomial $B$ 
together with (\ref{SUMPIF}) that
$\parallel \Phi_n \parallel^2 =\mathcal{O}(s_n^\beta)$ a.s. for all
$2\alpha^{-1}<\beta<1$. In addition, let $g_n=\Phi_n^tS_{n-1}^{-1}\Phi_n$
and $\delta_n=\text{tr}(S_{n-1}^{-1}-S_n^{-1})$. It is well-known that
$$(1-f_n)(1+g_n)=1$$ 
and $(\delta_n)$ tends to zero a.s.
Consequently, as $\parallel \pi_n \parallel^2=(1-f_n)(1+g_n)\parallel \pi_n \parallel^2$
and $1+g_n\leq 2 + \delta_n\parallel \Phi_n \parallel^2$, we can deduce from
(\ref{SUMPIF}) that
\begin{equation}  
\label{SPI}
\sum_{k=1}^n\parallel \pi_k \parallel^2=o(s_n^\beta\log s_n)
 \hspace{0.5cm} \text{a.s.}
\end{equation}
Therefore, we obtain from (\ref{CT}), (\ref{EQBASIS}) and (\ref{SPI}) that
\begin{equation}  
\label{SX}
\sum_{k=1}^n\parallel X_{k+1} \parallel^2=o(s_n^\beta\log s_n)+\mathcal{O}(n)
 \hspace{0.5cm} \text{a.s.}
\end{equation}
Furthermore, we infer from assumption $(\text{A}_1)$ that
\begin{equation}  
\label{SBU}
U_n=B^{-1}(R)A(R)X_{n+1}-B^{-1}(R)\veps_{n+1}
\end{equation}
which implies by (\ref{SX}) that
\begin{equation}  
\label{SU}
\sum_{k=1}^n\parallel U_{k} \parallel^2=o(s_n^\beta\log s_n)+\mathcal{O}(n)
\hspace{0.5cm} \text{a.s.}
\end{equation}
It remains to put together the two contributions (\ref{SX}) and (\ref{SU})
to deduce that $s_n=o(s_n)+\mathcal{O}(n)$ a.s. leading to
$s_n=\mathcal{O}(n)$ a.s. Hence, it follows from (\ref{SPI}) that
\begin{equation}  
\label{SP}
\sum_{k=1}^n\parallel \pi_k \parallel^2=o(n)
 \hspace{0.5cm} \text{a.s.}
\end{equation}
Consequently, we obtain from (\ref{CT}), (\ref{EQBASIS}) and (\ref{SP}) that
\begin{equation*}  
\lim_{n\rightarrow \infty} \frac{1}{n}\sum_{k=1}^n X_kX_k^t = \Gamma \hspace{0.5cm}\text{a.s.}
\vspace{-1ex}
\end{equation*}
and, for all $1\leq i \leq p-1$,
${\displaystyle \sum_{k=0}^n X_kX_{k-i}^t} = o(n)$ a.s.
which implies that
\begin{equation}
\label{CVGX}  
\lim_{n\rightarrow \infty} \frac{1}{n}\sum_{k=1}^n X_k^p(X_k^p)^t = L \hspace{0.5cm}\text{a.s.}
\end{equation}
where $L$ is given by (\ref{DEFL}). Furthermore, it follows from
relation (\ref{ARX}) together with (\ref{EQBASIS}) and 
assumption $(\text{A}_1)$ that, for all $n\geq 0$,
\begin{eqnarray*}  
U_n&=&B^{-1}(R)A(R)X_{n+1}-B^{-1}(R)\veps_{n+1}, \\
&=&B^{-1}(R)A(R)(\pi_n+x_{n+1})+B^{-1}(R)(A(R)-I_{d})\veps_{n+1}, \\
&=&V_n+W_{n+1}.
\end{eqnarray*}
Consequently, as $W_n=P(R)\veps_n$, we deduce from (\ref{CT}), (\ref{SP}) 
and the strong law of large numbers for martingales (see e.g. Theorem 4.3.16 of \cite{Duflo}) 
that, for all $1\leq i \leq q$,
\begin{equation*}  
\lim_{n\rightarrow \infty} \frac{1}{n}\sum_{k=1}^n U_kU_{k-i+1}^t = H_i \hspace{0.5cm}\text{a.s.}
\end{equation*}
which ensures that
\begin{equation}
\label{CVGU}  
\lim_{n\rightarrow \infty} \frac{1}{n}\sum_{k=1}^n U_k^q(U_k^q)^t = H \hspace{0.5cm}\text{a.s.}
\end{equation}
where $H$ is given by (\ref{DEFH}). Via the same lines, we also find that
\begin{equation}
\label{CVGXU}  
\lim_{n\rightarrow \infty} \frac{1}{n}\sum_{k=1}^n X_k^p(U_{k-1}^q)^t = K^t \hspace{0.5cm}\text{a.s.}
\end{equation}
Therefore, it follows from the conjunction of (\ref{CVGX}), 
(\ref{CVGU}) and (\ref{CVGXU}) that
\begin{equation}  
\label{CVGFIN}
\lim_{n\rightarrow \infty} \frac{S_{n}}{n} = \Lambda 
\hspace{0.5cm}\text{a.s.}
\end{equation}
where the limiting matrix $\Lambda$ is given by (\ref{DEFLAMBDA}). 
Hereafter, we recall that the $\mbox{ARX}_d(p,q)$ model is strongly
controllable. Thanks to Lemma \ref{MAINLEMMA}, the matrix $\Lambda$
is invertible and $\Lambda^{-1}$, given by (\ref{INVLAMBDA}), may be
explicitly calculated. This is the key point for the rest of the proof.
On the one hand, it follows from (\ref{CVGFIN}) that 
$n=\mathcal{O}(\lambda_{min}(S_n))$, $\parallel \Phi_n \parallel^2=o(n)$ a.s. 
which implies that $f_n$ tends to zero a.s. Hence, by  
(\ref{SUMPIF}), we find that
\begin{equation}  
\label{PIFIN}
\sum_{k=1}^n\parallel \pi_k \parallel^2=\mathcal{O}(\log n)
 \hspace{0.5cm} \text{a.s.}
\end{equation}
On the other hand, we obviously have from (\ref{EQBASIS})
\begin{equation}  
\label{COSTPI}
\parallel C_n - \Gamma_n \parallel =\mathcal{O}
\left(\frac{1}{n}\sum_{k=1}^n\parallel \pi_{k-1} \parallel^2\right)
\hspace{0.5cm} \text{a.s.}
\end{equation}
Consequently, we immediately obtain the tracking optimality (\ref{TH12})
from (\ref{PIFIN}) and (\ref{COSTPI}). Furthermore, by a well-known result
of Lai and Wei \cite{Lai} on the LS estimator, we also have
\begin{equation}
\label{RLS}
\parallel \wh{\theta}_{n+1} - \theta \parallel^{2}
= \mathcal{O} \left( \frac{\log \lambda_{max} S_{n}}{\lambda_{min} S_{n}}
\right) \hspace{0.5cm} \text{a.s.}
\end{equation}
Hence (\ref{TH14}) clearly follows from (\ref{CVGFIN}) and (\ref{RLS}).
Moreover, we also infer from Lemma 1 of Wei \cite{Wei} together with (\ref{CVGFIN}) that
\begin{equation}
\label{WEILS}
(\wh{\theta}_{n+1} - \theta)^tS_n(\wh{\theta}_{n+1} - \theta)
= o (\log n) \hspace{0.5cm} \text{a.s.}
\vspace{1ex}
\end{equation}
However, it follows from Theorem 4.3.16 part 4 of \cite{Duflo} that
\begin{equation}
\label{CVGDST}
\lim_{n\rightarrow \infty}
\frac{1}{\log d_n}\left(
(\wh{\theta}_{n+1} - \theta)^tS_n(\wh{\theta}_{n+1} - \theta)+
\sum_{k=0}^{n}(1-f_{k})\pi_{k}\pi_{k}^{t}\right)=
\Gamma
\hspace{0.5cm} \text{a.s.}
\end{equation}
where $d_n=det(S_n)$. In addition, if $\delta=d(p+q)$, we deduce from (\ref{CVGFIN}) that
\begin{equation}  
\label{CVGDET}
\lim_{n\rightarrow \infty} \frac{d_{n}}{n^\delta} = \det \Lambda 
\hspace{0.5cm}\text{a.s.}
\end{equation}
Finally, (\ref{WEILS}) together with (\ref{CVGDST}) and (\ref{CVGDET}) imply (\ref{TH13}), 
which achieves the proof of
Theorem \ref{ASPLS}. $\ \ \videbox$
\ \vspace{2ex}\par
\textit{Proof of Theorem \ref{ASPWLS}}. By Theorem 1 of \cite{Bercu1}, we have
\begin{equation}
\label{SUMPIFA}
\sum_{n=1}^{\infty} a_{n}(1-f_{n}(a))\parallel \pi_{n} \parallel^{2}
< + \infty \hspace{0.5cm} \text{a.s.}
\end{equation}
where $f_{n}(a)=a_{n}\Phi_{n}^{t}S_{n}^{-1}(a)\Phi_{n}$.
Then, as $a_{n}^{-1}=(\log(s_n))^{1+\gamma}$ with $\gamma>0$, 
we clearly have $a_{n}^{-1}=\mathcal{O}(s_n)$ a.s.
Hence, it follows from (\ref{SUMPIFA}) together with Kronecker's Lemma
given e.g. by Lemma 1.3.14 of \cite{Duflo} that
\begin{equation}
\label{SPIA}
\sum_{k=1}^{n}\parallel \pi_{k} \parallel^{2}
=o(s_{n}) \hspace{0.5cm} \text{a.s.}
\end{equation}
Therefore, we obtain from (\ref{CT}), (\ref{EQBASIS}) and (\ref{SPIA}) that
\begin{equation}  
\label{SXA}
\sum_{k=1}^n\parallel X_{k+1} \parallel^2=o(s_n)+\mathcal{O}(n)
 \hspace{0.5cm} \text{a.s.}
\end{equation}
In addition, we also deduce from assumption $(\text{A}_1)$ that
\begin{equation}  
\label{SUA}
\sum_{k=1}^n\parallel U_{k} \parallel^2=o(s_n)+\mathcal{O}(n)
\hspace{0.5cm} \text{a.s.}
\end{equation}
Consequently, we immediately infer from (\ref{SXA}) and (\ref{SUA})
that $s_n=o(s_n)+\mathcal{O}(n)$ so $s_n=\mathcal{O}(n)$ a.s.
Hence, (\ref{SPIA}) implies that
\begin{equation}  
\label{SPA}
\sum_{k=1}^n\parallel \pi_k \parallel^2=o(n)
 \hspace{0.5cm} \text{a.s.}
\end{equation}
Proceeding exactly as in Appendix A, we find from (\ref{SPA}) that
\begin{equation*}  
\lim_{n\rightarrow \infty} \frac{S_{n}}{n} = \Lambda 
\hspace{0.5cm}\text{a.s.}
\end{equation*}
Via an Abel transform, it ensures that
\begin{equation} 
\label{CVGFINA}
\lim_{n\rightarrow \infty} (\log n)^{1+\gamma} \frac{S_{n}(a)}{n} = 
\Lambda \hspace{0.5cm}\text{a.s.}
\end{equation}
We obviously have from (\ref{CVGFINA}) that $f_{n}(a)$ tends to zero a.s.
Consequently, we obtain from (\ref{SUMPIFA}) and Kronecker's Lemma
that
\begin{equation}
\label{PIFINA}
\sum_{k=1}^{n}\parallel \pi_{k} \parallel^{2}
=o((\log s_{n})^{1+\gamma}) \hspace{0.5cm} \text{a.s.}
\end{equation}
Then, (\ref{TH22}) clearly follows from (\ref{COSTPI}) and (\ref{PIFINA}).
Finally, by Theorem 1 of \cite{Bercu1} 
\begin{equation} 
\label{RWLS}
\parallel \wh{\theta}_{n+1} - \theta \parallel^{2}
=\mathcal{O} \left( \frac{1}{\lambda_{min} S_{n}(a)} \right)  
\hspace{0.5cm} \text{a.s.}
\end{equation}
Hence, we obtain (\ref{TH23}) from (\ref{CVGFINA}) and (\ref{RWLS}), which
completes the proof of Theorem \ref{ASPWLS}. $\ \ \videbox$

\section*{Appendix C.}

\renewcommand{\thesection}{\Alph{section}} \renewcommand{\theequation}{%
\thesection.\arabic{equation}} \setcounter{section}{3} 
\setcounter{equation}{0}

\textit{Proof of Theorem \ref{CLTLIL}}. First of all, it follows from (\ref%
{MOD}) together with (\ref{WLS}) that, for all $n\geq 1$, 
\begin{equation}  \label{DECWLS}
\widehat{\theta}_{n}-\theta=S_{n-1}^{-1}(a)M_{n}(a)
\end{equation}
where 
\begin{equation}  \label{DEFMART}
M_{n}(a)=\widehat{\theta}_0-\theta+\sum_{k=1}^n
a_{k-1}\Phi_{k-1}\varepsilon_k^t.
\end{equation}
We now make use of the CLT for multivariate martingales given e.g. by Lemma
C.1 of \cite{Bercu2}, see also \cite{Duflo}, \cite{Hall}. On the
one hand, for the LS algorithm, we clearly deduce (\ref{TH31}) from
convergence (\ref{TH11}) and decomposition (\ref{DECWLS}). On the other
hand, for the WLS algorithm, we also infer (\ref{TH31}) from convergence 
(\ref{TH21}) and (\ref{DECWLS}). Next, we make use of the LIL for
multivariate martingales given e.g. by Lemma C.2 of \cite{Bercu2},
see also \cite{Duflo}, \cite{Stout}. For the LS algorithm, since 
$(\varepsilon_n)$ has finite conditional moment of order $\alpha>2$, we
obtain from Chow's Lemma given e.g. by Corollary 2.8.5 of 
\cite{Stout} that, for all $2< \beta < \alpha$, 
\begin{equation}  \label{NORMEPS}
\sum_{k=1}^{n}\parallel \varepsilon_{k}\parallel^{\beta}=\mathcal{O}(n) 
\hspace{0.5cm}\text{a.s.}
\end{equation}
Consequently, as the reference trajectory $(x_n)$ satisfies (\ref{REGNORM}),
we deduce from (\ref{EQBASIS}) together with (\ref{PIFIN}) and (\ref{NORMEPS}) that 
\begin{equation}  \label{NORMX}
\sum_{k=1}^{n}\parallel X_{k}\parallel^{\beta}=\mathcal{O}(n) \hspace{0.5cm}
\text{a.s.}
\end{equation}
Furthermore, it follows from (\ref{SBU}) and (\ref{NORMX}) that 
\begin{equation}  \label{NORMU}
\sum_{k=1}^{n}\parallel U_{k}\parallel^{\beta}=\mathcal{O}(n) \hspace{0.5cm}
\text{a.s.}
\end{equation}
Hence, we clearly obtain from (\ref{NORMX}) and (\ref{NORMU}) that 
\begin{equation}  \label{NORMPHI}
\sum_{k=1}^{n}\parallel \Phi_{k}\parallel^{\beta}=\mathcal{O}(n) 
\hspace{0.5cm}\text{a.s.}
\end{equation}
Therefore, as $\beta >2$, (\ref{NORMPHI}) immediately implies that 
\begin{equation*}
\sum_{n=1}^{\infty}\left(\frac{\parallel \Phi_{n}\parallel} 
{\sqrt{n}}\right)^{\beta}<+\infty \hspace{0.5cm}\text{a.s.}
\end{equation*}
Finally, Lemma C.2 of \cite{Bercu2} together with convergence (\ref{TH11})
and (\ref{DECWLS}) lead to (\ref{TH32}). The proof for the WLS algorithm is
left to the reader because it follows essentially the same arguments than
the proof for the LS algorithm. It is only necessary to add the weighted
sequence $(a_n)$ and to make use of convergence (\ref{TH21}).$\ \ \videbox$


%

\end{document}